\documentstyle{amsppt}
\magnification =\magstep 1
\document

\topmatter
\title The distribution of vector-valued Rademacher series \endtitle
\author S.J. Dilworth and S.J. Montgomery-Smith \endauthor
\address Department of Mathematics, University of South Carolina,
Columbia, South Carolina 29208, U.S.A. \endaddress
\address Department of Mathematics, University of Missouri, Columbia,
Missouri 65211, U.S.A. \endaddress
\thanks The second author was supported in part by NSF DMS-9001796
\endthanks
\abstract Let $X=\sum \varepsilon_n x_n$ be a Rademacher series with
vector-valued coefficients. We obtain an approximate formula for the
 distribution
of the random variable $||X||$ in terms of its mean and a
certain quantity derived
from the
K-functional of
interpolation theory. Several applications of the formula are given.
\endabstract
\endtopmatter

\heading 1. Results \endheading
In \cite{6} the second-named author
calculated the distribution of a scalar Rademacher series $\sum
\varepsilon_n a_n$. The principal result of the present paper extends
the results of
\cite{6} to the case of a Rademacher series $\sum \varepsilon_n x_n$
with coefficients $(x_n)$ belonging to an arbitrary Banach space $E$.
Its proof relies on a deviation inequality for Rademacher series obtained by
Talagrand \cite{9}. A somewhat curious feature of the proof is that it
appears to exploit in a non-trivial way (see Lemma 2) the platitude that
every separable Banach space is isometric to a closed subspace of $\ell_\infty$.
The principal result is applied to yield a precise form
of the
Kahane-Khintchine
inequalities and to compute certain Orlicz norms for Rademacher
series.

First we recall some notation and terminology from interpolation theory
(see e.g. \cite{1}).
Let $(E_1,||.||_1)$ and $(E_2,||.||_2)$ be
two Banach spaces which are continuously embedded into some larger
topological vector space. For $t>0$, the K-functional $K(x,t;E_1,E_2)$ is
the norm on $E_1 + E_2$ defined by
$$ K(x,t;E_1,E_2) = \inf \lbrace ||x_1||_1 + t||x_2||_2: x=x_1+x_2,
\quad x_i \in E_i \rbrace. $$
For a sequence $(a_n) \in \ell_2$, we shall denote the K-functional
$K((a_n),t;\ell_1,\ell_2)$ by $K_{1,2}((a_n),t)$ for short. For
$1\le p<\infty$, a sequence $(x_n)$ in a Banach space $(E,||.||)$ is said to
be weakly-$\ell_p$ if the scalar sequence $(x^*(x_n))$ belongs to $\ell_p$
for every $x^* \in E^*$. The collection of all weakly-$\ell_p$ sequences is a
Banach space, denoted $\ell^w_p(E)$, with the norm given by
$l^w_p((x_n))=\sup_{||x^*||\le 1} ||(x^*(x_n))||_p$ (where
$||(a_n)||_p = (\sum |a_n|^p)^{1/p}$). If $(x_n)\in
\ell^w_2(E)$, we make the following definition:
$$K^w_{1,2}((x_n),t)= \sup_{||x^*||\le 1} K_{1,2}((x^*(x_n)),t).$$
Observe that $K^w_{1,2}((x_n),t)$ is a continuous increasing function of $t$.
In fact, it is a Lipschitz function with Lipschitz constant at most $\ell
^w_2((x_n))$.

Next we set up some function space notation. Let $(\Omega,\Sigma,P)$ be
a probability space. A Rademacher (or Bernoulli) sequence $(\varepsilon_n)$ is
a sequence of independent identically distributed random variables such that
$P(\varepsilon_n = 1)= P(\varepsilon_n = -1) = \frac12$. For a random variable
$Y$ defined on $\Omega$, its decreasing rearrangement, $Y^*$, is the function on
$[0,1]$ defined by $Y^*(t)= \inf \lbrace s>0 : P(|Y|>s) \le t \rbrace$.
For $0<p<\infty$, the weak-$L_p$ norm of $Y$, denoted $||Y||_{p,\infty}$, is
 given by
$||Y||_{p,\infty} = \sup_{0<t<1}t^{\frac 1p}Y^*(t)$. As usual, $||Y||_p$
denotes $(\Bbb E |Y|^p)^{1/p}$.
Let $\Psi$ be an Orlicz function on $[0,\infty)$. The Orlicz norm, $||Y||_\Psi$
, is given by $||Y||_\Psi = \inf \lbrace c>0: \Bbb E \Psi(|Y|/c) \le 1 \rbrace.$
 We shall be particularly interested in the Orlicz functions
 $\Psi_q(t)=e^{t^q}-1$ for $2<q<\infty$. The weak-$\ell_p$ norm of the scalar
 sequence $(a_n)$ is
defined by $||(a_n)||_{p,\infty}= \sup n^{\frac 1p} a^*_n$, where $(a_n^*)$
is the decreasing rearrangement of $(|a_n|)$.

Finally, we shall write $A\approx B$ to mean that there is a constant $C>0$ such
 that $\frac1C A \le B \le CA$.
We shall try to indicate in each case whether the implied constant is absolute
or whether it depends on some parameter, typically $p \in [1,\infty)$, entering
into the expressions for $A$ and $B$.

Now we can state the principal result of the paper.
\proclaim{\smc Main Theorem} Let $X=\sum \varepsilon_n x_n$ be an almost surely
convergent Rademacher series in a Banach space $E$. Then, for $t>0$, we have
$$P(||X|| > 2\Bbb E ||X|| + 6 K^w_{1,2}((x_n),t)) \le 4e^{-t^2/8},
\tag1$$
and, for some absolute constant $c$, we have
$$P\left(||X||>\frac12 \Bbb E ||X|| + cK^w_{1,2}((x_n),t)
\right) \ge c e^{-t^2/c}.
\tag2$$ \endproclaim
The proof of the Main Theorem will be deferred until the end of the paper in
order to proceed at once with the applications.
\proclaim{Corollary 1} Let $X=\sum \varepsilon_n x_n$ be an almost surely
convergent Rademacher series in a Banach space. Then, for $0<t\le \frac{1}{10}$,
 we have
$$S^*(t) \approx \Bbb E ||X|| + K^w_{1,2}((x_n),\sqrt{\log (1/t)}), \tag3$$
where $S$ denotes the real random variable $||X||$. The implied constant is
absolute. \endproclaim
\demo{Proof} (1) and (2) give rise to the inequalities $S^*(4e^{-t^2/8})
\le 2\Bbb E ||X|| + 6K^w_{1,2}((x_n),t)$ and $S^*(ce^{-t^2/c}) \ge \frac12
\Bbb E ||X|| + cK^w_{1,2}((x_n),t)$, respectively, whence (3) follows for
all sufficiently small $t$ by an appropriate change of variable. To see that
the lower estimate implicit in (3) is valid in the whole range $0<t< \frac{1}
{10}$, we recall from \cite{2} that $\Bbb E ||X||^2 \le 9 \Bbb E^2||X||$.
Hence, by the Paley-Zygmund inequality
(see e.g. \cite{4,p.8}), for $0<\lambda<1$, we have
$$\align P(||X|| > \lambda \Bbb E ||X||) &\ge (1-\lambda)^2 \frac{\Bbb E^2X}
{\Bbb E X^2} \\
&\ge \frac19 (1-\lambda)^2, \endalign$$
whence $P(||X|| > (1- \frac{3}{\sqrt{10}})\Bbb E ||X||) \ge \frac{1}{10},$
which easily implies (3). \qed \enddemo
In \cite{4} Kahane proved that if $P(||X|| > t) = \alpha$, where $X$ is a
Rademacher series in a Banach space, then $P(||X|| > 2t) \le 4\alpha^2.$
By iteration this implies $P(||X|| > st)\le \frac14 (4\alpha)^s$ for $s=
2^n$.  According to our next corollary the exponent $s$ in the
latter result may be
improved to be a certain multiple of $s^2$.
\proclaim{Corollary 2} Let $X = \sum \varepsilon_n x_n$ be an almost surely
convergent Rademacher series in a Banach space. Then, for $t>0$ and
$s \ge 1$, we have
$$P(||X|| > st) \le \left(\frac{1}{c_1}P(||X||>t)\right)^{c_1s^2}$$
for some absolute constant $c_1$. \endproclaim
\demo{Proof} By choosing $c_1<c$, where $c$ is the constant which appears
in (2), the result becomes trivial whenever $P(||X||>t)\ge c.$ Hence
we may assume that $P(||X||>t)<c.$ Choose $\alpha > 0$ such that
$P(||X||>t)= c e^{-\alpha^2/c}.$ Then (2) gives $t\ge \frac12 \Bbb E ||X|| +
cK^w_{1,2}((x_n),\alpha).$ Thus,
$$ \align st &\ge \frac
{s}{2}  \Bbb E ||X|| + sc K^w_{1,2}((x_n),\alpha) \\
&\ge 2 \Bbb E ||X|| + K^w_{1,2}((x_n),cs\alpha) \endalign $$
provided $s\ge \max(4,1/c)$. Now (1) gives
$$\align P(||X||>st) &\le 4e^{-(cs\alpha)^2/8} \\
&= 4\left(\frac1c(ce^{-\alpha^2/c})\right)^{c^3s^2/8} \\
&= 4\left(\frac1c(P(||X|| > t)\right)^{c^3s^2/8}, \endalign $$
which gives the result. \qed \enddemo
Our next corollary,
which is the vector-valued version of a recent result of
Hitczenko \cite{3}, is a rather precise form of the Kahane-Khintchine
 inequalities.
\proclaim{Corollary 3} Let $X = \sum \varepsilon_n x_n$ be a Rademacher series
in a Banach space. Then, for $1\le p <\infty$, we have
$$(\Bbb E ||X||^p)^{1/p} \approx \Bbb E ||X|| + K^w_{1,2}((x_n),\sqrt p).$$
The implied constant is absolute. \endproclaim
\demo{Proof} We may assume that $p\ge2$. It follows from a result of
Borell \cite{2} that \linebreak
$(\Bbb E ||X||^{2p})^{1/2p} \le \sqrt3(\Bbb E ||X||^p)^{1/p}.$
Since $\frac12 ||Y||_p \le ||Y||_{2p,\infty} \le ||Y||_{2p}$ for every random
variable $Y$ (as is easily verified),
it follows (letting $S$ denote the random variable $||X||$)
that $\frac12 ||S||_p \le ||S||_{2p,\infty} \le \sqrt3 ||S||_p$.
So it  suffices to prove that $||S||_{p,\infty} \approx \Bbb E S +
K^w_{1,2}((x_n),\sqrt p)$ to obtain the desired conclusion.
By Corollary 1, we have
$$ \align
||S||_{p,\infty} &\approx \Bbb E S + \sup_{0<t<1}t^{1/p}K^w_{1,2}((x_n),
\sqrt{\log(1/t)}) \\
&=\Bbb E S + \sup_{0<t<1}\left\{
 t^{1/p}\sup_{||x^*||\le 1} K_{1,2}((x^*(x_n)),\sqrt{\log (1/t)
}) \right\} \\
&= \Bbb E S + \sup_{||x^*||\le 1} \left\{ \sup_{0<t<1} t^{1/p} K_{1,2}
((x^*(x_n)),\sqrt{\log (1/t)}) \right\}. \endalign $$
To evaluate the expression in brackets we shall make use once more (see
Corollary 2) of the elementary inequality
$K_{1,2}((a_n),s) \le \max (1,s/t) K_{1,2}((a_n),t).$ Thus,
$$\align \sup_{0<t\le e^{-p}}t^{1/p} K_{1,2}((x^*(x_n)),\sqrt{\log(1/t)})
&\le \left(\sup_{0<t\le e^{-p}} t^{1/p} \sqrt{\dfrac {\log(1/t)}{p}}\right)
K_{1,2}((x^*(x_n)),\sqrt p) \\
&= e^{-1} K_{1,2}((x^*(x_n)),\sqrt p). \endalign $$
Moreover,
$$\sup_{e^{-p}<t<1} t^{1/p} K_{1,2}((x^*(x_n)),\sqrt{\log (1/t)})
\le K_{1,2}((x^*(x_n)),\sqrt p).$$
Finally, we obtain
$$\frac1e K^w_{1,2}((x_n),\sqrt p) \le \sup_{||x^*|| \le 1} \left\{
\sup_{0<t<1} K_{1,2}((x^*(x_n)),\sqrt{\log (1/t)}) \right\} \le
K_{1,2}^w((x_n),\sqrt p),$$
which gives the desired result.
\qed \enddemo
Our final application is to the calculation of the Orlicz norms $||S||_{\psi
_q}$ for $2<q<\infty.$ The proof will use the scalar version of the result,
which
was obtained by Rodin and Semyonov \cite{8} (see also \cite{7}). (Recall that by
a result of Kwapien, \cite{5}, $||S||_{\psi_q} \approx
\Bbb E ||X||$ in the range $0<q\le2$.)
\proclaim{Corollary 4} Let $X=\sum \varepsilon_n x_n$ be an almost surely
convergent Rademacher series in a Banach space. Then, for $2<q<\infty$, we
have
$$||S||_{\psi_q} \approx \Bbb E ||X|| + \sup_{||x^*|| \le 1} ||(x^*(x_n))||_
{p,\infty},$$
where $\frac1p + \frac1q = 1$ and $S$ denotes $||X||$. The implied constant
depends only on $q$. \endproclaim
\demo{Proof} It is easily verified that $||f||_{\psi_q} \approx \sup_{0<t<1}
(\log(1/t))^{-1/q} f^*(t). $ Hence, by Corollary~1, we have $$\align
||S||_{\psi_q} &\approx \Bbb E ||X|| + \sup_{0<t<1}(\log(1/t))^{-1/q}K^w_{1,2}
((x_n),t) \\
&\approx \Bbb E ||X|| + \sup_{0<t<1}\left\{ (\log(1/t))^{-1/q} \sup_{||x^*||\le
 1}
K_{1,2}((x^*(x_n)),t) \right\} \\
&\approx \Bbb E ||X|| + \sup_{||x^*||\le 1} \left\{ \sup_{0<t<1}
 (\log(1/t))^{-1/q}
K_{1,2}((x^*(x_n)),t) \right\} \\
&\approx \Bbb E ||X|| + \sup_{||x^*|| \le 1}||\sum \varepsilon_n
 x^*(x_n)||_{\psi_q}\\
&\approx \Bbb E ||X|| + \sup_{||x^*|| \le 1} ||(x^*(x_n))||_{p,\infty},
 \endalign $$
where the last line follows from the result of Rodin and Semyonov.
\qed \enddemo

\heading 2. Proof of main result \endheading
The principal ingredient in the proof of the Main Theorem is the following
deviation inequality of Talagrand \cite{9}.
\proclaim{Theorem A} Let
$X=\sum_{n=1}^N \varepsilon_n x_n$ be a finite Rademacher series in a Banach
space and let $M$ be a median of $||X||$. Then, for $t>0$, we have
$$ P\left(|||\sum_{n=1}^N \varepsilon_n x_n||-M|>t\right)
\le 4 e^{-t^2/8\sigma^2},
$$ where $\sigma = \ell^w_2((x_n)_{n=1}^N).$ \endproclaim
\proclaim{Lemma 1} Let $X=\sum_{n=1}^N \varepsilon_n x_n$ be a finite Rademacher
 series in a Banach space $E$. Then, for $t>0$, we have
$$P(||X||>2\Bbb E||X|| + 3K((x_n)_{n=1}^N,t;\ell^w_1(E),\ell^w_2(E))) \le
4e^{-t^2/8}.$$ \endproclaim
\demo{Proof} It follows from Theorem A that for all $y_1, \dots, y_N$ in $E$,
we have
$$P(||\sum\varepsilon_n y_n||>2\Bbb E||\sum \varepsilon_n y_n||
+ t \ell^w_2((y_n)))\le 4e^{-t^2/8}. \tag4$$
On the other hand, since $\max ||\sum \varepsilon_n y_n|| = \ell^w_1((
y_n))$, we have the trivial estimate
$$P(||\sum \varepsilon_n y_n|| > \ell^w_1((y_n))) = 0. \tag5$$
Let $x_n= x^{(1)}_n + x^{(2)}_n$ for $1\le n\le N$, let $X^{(1)}= \sum
 \varepsilon_n x^{(1)}_n$, and let $X^{(2)}= \sum \varepsilon_n x^{(2)}_n.$ Then
$$ \align
\ell^w_1((x^{(1)}_n)) + t \ell^w_2((x^{(2)}_n)) + 2\Bbb E ||X^{(2)}|| &\le
\ell^w_1((x^{(1)}_n)) + t \ell^w_2((x^{(2)}_n)) + 2\Bbb E
||X^{(1)}|| + 2\Bbb E ||X|| \\
&\le 3\ell^w_1((x^{(1)}_n)) + t\ell^w_2((x^{(2)}_n)) + 2\Bbb E ||X|| \\
&\le 2 \Bbb E ||X|| + 3(\ell^w_1((x_n^{(1)})) + t\ell^w_2((x^{(2)}_n))).
\endalign $$
Let $Q$ denote $2\Bbb E ||X|| + 3 (\ell^w_1((x^{(1)}_n)) + t
 \ell^w_2((x_n^{(2)})))$. Then, by (4) and (5) and by the above inequality, we
 have
$$\align P(||X|| > Q) &\le P(||X^{(1)}||+ ||X^{(2)}|| > \ell^w_1((x_n^{(1)}))
+ t \ell^w_2((x^{(2)}_n))
+ 2 \Bbb E ||X^{(2)}||) \\
&\le P(||X^{(1)}|| > \ell^w_1((x_n^{(1)})))
+ P(||X^{(2)}|| > 2\Bbb E ||X^{(2)}||
+ t \ell^w_2((x^{(2)}_n))) \\
&< 0 + 4e^{-t^2/8} \endalign$$
The desired conclusion now follows from the definition of the K-functional.
\qed
\enddemo
\proclaim{Lemma 2} Let $x_1, \dots, x_N$ be elements of the Banach space
$\ell_\infty.$ Then
$$K((x_n)_{n=1}^N,t; \ell^w_1(\ell_\infty),\ell^w_2(\ell_\infty)) \le
2 K^w_{1,2}((x_n)_{n=1}^N,t).$$ \endproclaim
\demo{Proof}
For $1\le n\le N$, let $x_n = (x_{n,j})_{j=1}^\infty \in \ell_\infty.$
A simple convexity argument gives
$$||(x_n)||_{\ell^w_p(\ell_\infty)} = \sup_{1 \le j \le \infty}
\left(\sum_{n=1}^N |x_{n,j}|
^p\right)^{(1/p)}.$$
It follows that the mapping $\phi$ which associates an element $(y_n)_{n=1}
^\infty \in \ell^w_p(\ell_\infty)$ with the element in $\ell_\infty(\ell_p)$
whose $jth$ coordinate equals $(y_{n,j})_{n=1}^\infty$ is an isometry.
Hence \linebreak
$K((x_n),t;\ell^w_1,\ell^w_2)= K(\phi((x_n)),t;\ell_\infty(\ell_1),
\ell_\infty(\ell_2)).$  Let $(y_n)_{n=1}^\infty \in \ell_\infty(\ell_2)$
and let $\varepsilon >0$. For each $n$ there exists a splitting
$y_{n} = z^{(1)}_{n} + z^{(2)}_{n}$ such that
$$||(z^{(1)}_{n,j})_{j=1}^\infty||_1 + t||(z^{(2)}_{n,j})_{j=1}^\infty||_2
\le K_{1,2}((y_{n,j})_{j=1}^\infty,t) + \varepsilon.$$
It follows that $$ \align
||(z^{(1)}_n)||_{\ell_\infty(\ell_1)} + t||(z^{(2)}_n)||_{\ell_\infty(\ell
_2)} &= \sup_{1 \le n < \infty} ||(z^{(1)}_{n,j})_{j=1}^\infty||_1
+ t\sup_{1 \le n < \infty}||(z^{(2)}_{n,j})_{j=1}^\infty||_2 \\
&\le 2\sup_{1 \le n < \infty} K_{1,2}((y_{n,j})_{j=1}^\infty,t) + 2\varepsilon
\\ &\le 2 K^w_{1,2}((y_n),t) + 2 \varepsilon. \endalign$$
Since $\varepsilon$ is arbitrary, the result now follows from the definition of
 the
K-functional. \qed
\enddemo
\demo {Proof of Main Theorem} First we prove (1) for a finite Rademacher
series $\sum_{n=1}^N \varepsilon_n x_n.$ Since every separable Banach space
embeds isometrically into $\ell_\infty$, we may assume that $E$ is a closed
subspace of $\ell_\infty$. Recall that $K^w_{1,2}((x_n),t)$ was defined as
$\sup_{||x^*|| \le 1}K_{1,2}((x^*(x_n)),t).$ By the Hahn-Banach Theorem, the
supremum is the same whether it is taken over elements of $E^*$ or
over elements of
$\ell_\infty^*$. Hence (1) follows by combining Lemmas 1 and 2. The result
for an infinite series follows from the result for $\sum_{n=1}^N \varepsilon_n
x_n$ by taking the limit as $N\rightarrow \infty.$ To prove (2), we use the
result from \cite{6} that there exists an absolute constant $d$ such that
$P(\sum \varepsilon_n a_n > d K_{1,2}((a_n),t)) \ge d e^{-t^2/d}$ for every
sequence $(a_n) \in \ell_2$. Hence
$$\align
P\left(||\sum \varepsilon_n x_n|| >  \frac{d}{2} K^w_{1,2}((x_n),t)\right) &\ge
\inf_{||x^*|| \le 1} P(||\sum \varepsilon_n x_n|| > d K_{1,2}((x^*(x_n)),t))\\
&\ge \inf_{||x^*|| \le 1} P(\sum \varepsilon_n x^*(x_n) > d
 K_{1,2}((x^*(x_n)),t))\\
&\ge de^{-t^2/d}. \endalign $$
The Paley-Zygmund inequality now gives
$$ \align
P\Biggl(||X||> \frac12 \Bbb E ||X||
&+ \frac{d}{6} K^w_{1,2}((x_n),t) \Biggr) \\
&\ge
\min \left( P \left(||X||>\frac34 \Bbb E ||X|| \right),
P\left(||X|| > \frac {d}{2} K^w_{1,2}((x_n),t)\right) \right) \\
&\ge \min \left( \frac{1}{144}, de^{-t^2/d}. \right). \qed \endalign $$
\enddemo

\enddocument